# Optimizing Daily Fantasy Baseball Lineups: A Linear Programming Approach for Enhanced Accuracy

Max Grody[1], Sandeep Bansal[1], Huthaifa I. Ashqar[1,2]

[1] University of Maryland, Baltimore County

[2] Arab American University

**Abstract**

Daily fantasy baseball has shortened the life cycle of an entire fantasy season into a single day. As of today, it has become familiar with more than 10 million people around the world who participate in online fantasy. As daily fantasy continues to grow, the importance of selecting a winning lineup becomes more valuable. The purpose of this paper is to determine how accurate FanDuel's current daily fantasy strategy of optimizing daily lineups are and utilize python and linear programming to build a lineup optimizer for daily fantasy sports with the goal of proposing a more accurate model to assist daily fantasy participants select a winning lineup.

**Introduction**

Daily fantasy baseball has revolutionized the way people participate in online fantasy baseball leagues. With the compression of an entire season into a single day, sites like FanDual have provided fans with the unique experience of selecting a new lineup every day. A lineup also referred to as an entry holds constraints on players used to build any given team. FanDual, a daily fantasy sports company sets a virtual price for each player. These prices differ daily and are influenced by past and predicted performances.

An added challenge includes a maximum number of each type of player allowed per entry. For example a daily fantasy lineup consists of a maximum of 9 players made up of one pitcher, catcher/first, second, and third baseman, three outfielders, a shortstop and a utility position which allows the participant to select a player of any position to be a part of their entry. Although these constraints increase the difficulty of selecting a winning lineup, the largest challenge remains: How does one select and optimize their entry if the performance of the players is unknown? In other words, how does one know the player he or she selects will perform well in order to increase their chances of winning?

Along with providing participants with a platform to select and play their entries Fandual also provides a predictive model of how well the entries will perform. A participant may utilize Fandual's predictive model to assist selecting entries that have a high probability of winning however, due to the lack of consistency and accuracy of Fandual fantasy platform participants typically utilize multiple resources in order to select the best entry. Ultimately participants must distinguish the trade-off between probability of entries and their correlation in order to select a good set of entries.

Guaranteed Prize Pool (GPP) tournaments is regarded as the most popular contest in Daily Fantasy Sports. FanDual guaranteed prizes even if the contest isn't completely full. In order to be victorious a participant is recommended to utilize entries that differ from other participants. These contests occur every day, allowing multiple opportunities to test and evaluate different strategies. For these reasons, daily fantasy sports provide a great opportunity for participants to test their strategy to picking a winning lineup.

**Research Question**

The underlying question we would like to answer is: How accurate is Fandual in regard to their accuracy for daily fantasy? More importantly are we able to utilize python and linear programming to assist participants in building their own lineup? Linear programming is an attempt to find a maximum or minimum solution to a problem based on a combination of variables subject to a given set of constraints. In order to explore and analyze we first have to define variables as players to include in a lineup. These variables would include a few defining features such as their position, salary, and projection. Next we would add in our constraints. Some basic constraints include things such as ensuring the lineup stays under the total salary cap, and positional requirements. We then define the objective which is to maximize a lineup's projected points.

**Literature Review**

Since daily fantasy sports is a relatively new industry, there has been minimal research exploring strategies to increase winning probability. One of the first papers published on such a topic is *Optimizing daily fantasy sports contests through stochastic integer programming* by Sarah Newell. In this paper, author Sarah Newell explores a stochastic integer programming approach to optimizing the expected payout of any given lineup in daily fantasy football. Specifically, she determines what scores a lineup needs to hit to be expected to win a payout. For example, she estimates that for one of the million dollar contests, in order to win the million dollars, a lineup would need to score at least 253.7 fantasy points, and in order to win the second place prize of $100,000, a lineup would need to score at least 250.57 fantasy points (cite this, page 26). She then calculates the expected payout of a lineup by taking the sum of the probability that the lineup hits each benchmark score times the increase in prize value from the

previous prize. Her model is then designed to optimize lineups to maximize expected payout. Her model falls short in that the single lineups it produces do not perform well. This is in part due to the limitations of her input data, which were simply the average of a player's scores from all previous games. While this is a much simpler projection system than our projection data, it is likely less accurate at predicting future performance. Additionally, single entry lineups are unlikely to win tournaments given the hundreds of thousands of entries each daily tournament typically receives.

In the paper *Picking Winners in Daily Fantasy Sports Using Integer Programming*, the authors build off of the integer programming approach of Sarah Newell by introducing the concepts of "stacking" lineups and creating diverse lineup portfolios for daily fantasy hockey. In Sarah Newell's paper she makes the assumption that a player's fantasy points are normally distributed and independent of other players. The authors of *Picking Winners* challenge this assumption and provide evidence that, at least for hockey, there is a strong positive correlation for fantasy points between players on the same team. This discovery is implemented in their integer programming model such that the model will produce a portfolio of optimal lineups in which projected fantasy points are maximized and include players from the same team. The authors' goal is to increase the variance of each lineup as much as possible in order to increase the probability that one of the lineups will score enough to win a top tier prize. Additionally, constraints were added to the model to diversify the lineup portfolio, ensuring that each lineup is different from each other, so as to capture more of the probability space and boost the probability that at least one lineup in the portfolio wins. Our approach to optimizing lineups is very similar to this one, with the exception that we are testing this approach for baseball instead of hockey. Additionally, further research will be pursued to improve the model. Possibilities include changing the input data to increase each lineup's variance and including data that forecasts what percent of lineups each player is on in a given tournament (this can help identify players who have both a high chance of scoring very well and who are not owned by many other lineups in the contest).

In *How to Play Strategically in Fantasy Sports (and Win)*, the authors try to improve on the success of the model from *Picking Winners*.  The authors of this paper are interested in winning both GPP tournaments and other types of daily fantasy tournaments that do not have the same type of top-heavy payout structure.  Additionally, the authors develop a way to reduce the optimization problem to a set of binary quadratic programs while also using Dirichlet regression to model player selections made by other contestants in order to improve the success of their optimized lineups against the competition.

Of the limited number of publications on fantasy sports strategy, it is apparent that optimization models pertaining to concepts of linear and quadratic programming are a popular strategy.  Given this information, creating a similar type of optimization model for daily fantasy baseball makes intuitive sense.

**Dataset**

The dataset was retrieved from a third party website named 'Sabersim' which is a website that supports Fandual daily fantasy and is able to provide participants with individual player statistics and projections. A few challenges that came up were how to obtain data that would accurately describe a player projected and accurate performance in order to assess Fandual's method, and to further explore and test our method for comparison. Upon review we decided to pick a timeframe in which players would at their best and would depict a precise description of how each player performed. Ultimately, we would have liked to select an entire season's worth of data but due to the time constraint we selected June 1- June 11 from the 2019 season.

While acquiring the data, we faced another challenge in which multiple players' actual score was not recorded and were left with a blank. After further investigation we concluded that these players were of the lower tiered in which they did not enter the game and were left with a null value as their actual score. Ultimately, we decided to drop those players entirely since those were players participants would be highly unlikely to choose to represent their entry in daily fantasy. Our complete dataset was represented by 408 total players, and four different columns including the player's projection score, actual score, the price at which they cost, and the difference between the projection and actual score.

After finalizing our dataset and completing the descriptive analysis portion we concluded that the average projection for each player is approximately 8.41 points, and the average actual score was approximately 9.55 points. The average price per player was about $3903. One observation to be made is that the actual minimum was -15. This depicts that a player performed poorly and ended up with a negative output. See figure 1

|       | SS Projection | Actual      | Price        | difference  |
|-------|---------------|-------------|--------------|-------------|
| count | 408.000000    | 408.000000  | 408.000000   | 408.000000  |
| mean  | 8.418504      | 9.558333    | 3903.186275  | 110.476588  |
| std   | 7.303810      | 11.746994   | 1681.132368  | 233.517723  |
| min   | 0.000000      | -15.000000  | 2000.000000  | 0.000000    |
| 25%   | 2.133050      | 3.000000    | 2600.000000  | 5.591419    |
| 50%   | 9.000660      | 6.000000    | 3300.000000  | 37.254126   |
| 75%   | 11.018675     | 12.775000   | 5500.000000  | 93.971630   |
| max   | 45.476800     | 73.000000   | 11500.000000 | 1914.893840 |

Figure 1.

After our initial analysis we calculated the mean squared error in order to determine the difference of the actual. The calculation was equal to 10.510 points. This number reflects how much on average the actual differed from the projection. We then wanted to visually describe how each player's actual differed from the projected, in figure 2 the reader is able to visualize each player and how they are are projected in comparison to how they actually performed. The arrow specifically is illustrating how FanDual projected a player to score roughly 12 points, but actually scored roughly 3 points, netting a difference of approximately 9 points.

Figure 2:

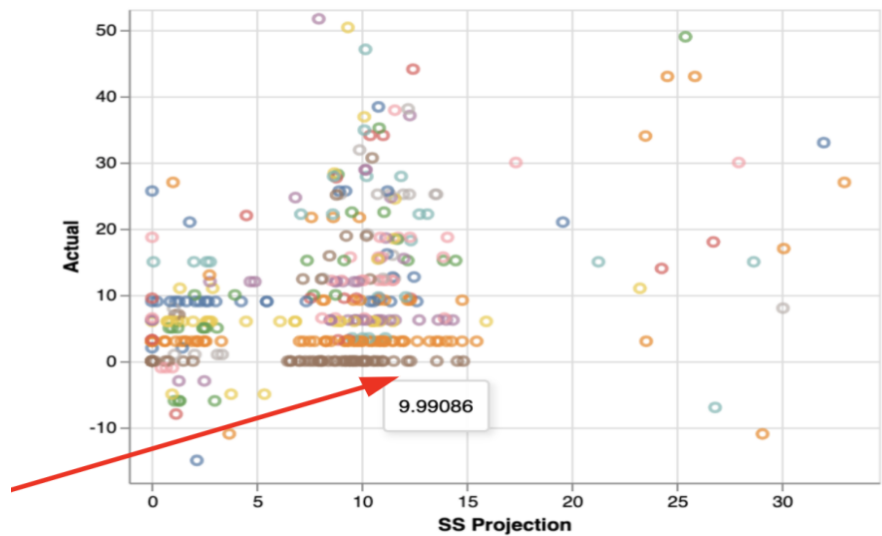

Next we wanted to illustrate how accurate Fandual was over the time period in which we selected. Ultimately we wanted to know: Are Fanduel projections reliable? In order to do we calculated the R-mean-squared which equated to 0.19. This number ultimately confirms our null hypothesis over the period of time we selected. Since the R-mean-squared was significantly less that 1, it is quantitatively significant that Fandual's projection mode is not reliable for those who participate in daily fantasy baseball.

**Methods**

Our problem is how to create a portfolio of lineups that will increase a contestant's probability to win a daily fantasy baseball contest on any given day or night. The team that scores the most fantasy points wins the contest, so we know that in order to win, we need to maximize the amount of fantasy points our lineups score. This is an optimization problem, as our goal is to maximize fantasy points. Since we are working with projection data, we will be trying to maximize the amount of projected fantasy points for each lineup we create. In the Picking Winners article, the authors decide to use an integer programming algorithm, and this is the strategy we implemented as well. Integer programming is ideal for this specific type of problem. One of the most well-known problems of combinatorial optimization is "the knapsack problem." In the knapsack problem, one is challenged with packing a bag with items. Each item has a specific value of utility to the user and a weight. The goal of the problem is to fit as much value into the bag without exceeding the bags weight limit. This problem sounds very similar to daily fantasy baseball. Players offer a specific value (projected fantasy points), and in this case their "weight" would be their individual prices. A popular solution to the knapsack problem is integer programming, so that makes it an ideal option for daily fantasy baseball, given the similarities.

Integer programming problems have several components: an objective function, constraints, and decision variables. The objective function is the function that maximizes or minimizes. For daily fantasy baseball, the objective function is to maximize projected fantasy points. Alternatively, the constraints limit the field of possible combinations. In daily fantasy baseball, lineups are constrained by how much they cost (a lineup cannot exceed the cost of $35000), the positions that players play, and the number of players in the lineup. Finally, the

decision variables are the variables that the algorithm decides to take or not. For fantasy baseball, the decision variables indicate whether a player is included in the lineup or not. An integer programming approach is more suited for daily fantasy baseball than a simple linear programming approach because our decision variables must be integers. It is not possible in fantasy sports to pay for a part of a player's fantasy points, either the contestant pays the price to put the player in their lineup, or they do not put the player in the lineup, those are the only possibilities. To reflect that, the decision variables of the model were made binary. A player is assigned a 1 if they are to be included in the lineup and a 0 if they are not to be included.

      The integer programming model was designed using Python and the PuLP library, which offers linear programming capabilities for Python. PuLP offers several options for solvers (some of which require subscriptions), including GUROBI and CPLEX, but the free default COIN solver was effective and fast enough for our purposes. Running this model successfully produces a single lineup optimized to maximize projected fantasy points given budget and positional constraints.

      Unfortunately, given the top-heavy payout structure of GPP tournaments (in which the largest prizes are offered to the lineups that score the most points; all other lineups win little to no money) and the fact that daily fantasy contests typically allow up to 150 entries into a single contest, it makes very little sense to only enter one lineup. Contests typically have hundreds of thousands of entries, and despite our single lineup being optimized as best it can given the limited information we have, it is highly unlikely that our lineup will be the one to score the most points out of the hundreds of thousands of other lineups. In order to make our model more practical, we adjusted our model to run iteratively to produce multiple optimal lineups. The more lineups a user enters into the contest, the higher the probability that one of the user's

lineups wins. As long as one of the user's lineups wins, positive profit margins will be achieved, since the top prizes range from $10,000-$50,000 and entry fees for each lineup are usually about $5. Specifically for the model, a constraint was added to make the next lineup iteration marginally smaller than the previous iteration, therefore ensuring that a new combination of players is produced each iteration.

An issue that we encountered was that all of the lineups produced by the model looked very similar to each other. Each new lineup would only change 1 or 2 players from the previous lineup. Strategically this is problematic. If all of the lineups are close to identical or similar in composition, this limits the probability that at least one of those lineups will succeed. If one lineup fails, then the rest are also likely to fail since they are all built of essentially the same players. To counter this problem, we adjusted the model to remove players as options once they had appeared on a set number of lineups. For example, if the model is run to produce 100 lineups and the user does not want a single player to appear on more than 25 lineups, once a player appears on their 25th lineup, they will then no longer be considered for the remaining lineups to be produced. This adjustment allows a lineup portfolio to be more diverse and thus increase the portfolio's probability of including at least one successful lineup.

**Analysis and Results**

In order to have access to past daily fantasy baseball contests, one must have entered into that contest. Since as of this writing there are no contests being offered, the model cannot yet properly be assessed of its success against other contestants (past or present). However, there is still insight to be gained from running the model on past 2019 projection data. Over a 30 day span, the model's lineup with the most projected fantasy points was projected to score on average 107.3 fantasy points fewer than the amount of points actually scored for the most optimal lineup for that day. Specifically, the highest projected lineups over those 30 days were projected to score on average 144.6 points and the actual optimal lineup for each day on average scored 251.9 points. That is a large discrepancy. The discrepancy is more a fault of the input data than it is a negative reflection of the integer programming approach. The actual optimal lineups consist almost exclusively of players that significantly outperformed their projections, and since our model only produces lineups based on projections (not the scores players would achieve if they exceeded their projections), lineups produced by any optimizer will never have a projected score higher than the actual optimized lineup.

It will be important to test the model once daily fantasy baseball contests return in order to truly evaluate how it performs against other teams. Additionally, there are many next steps to be taken to further increase the probability of winning a tournament.

**Conclusion**

Ultimately, we were able to create an integer programming model capable of producing multiple lineup entries optimized for projected fantasy points in a daily fantasy MLB contest. The integer programming approach pursued for daily fantasy hockey by the authors of "Picking Winners" is in fact generalizable to a sport like baseball.

While the current iteration of this model would provide an edge over amateur daily fantasy players, there are still many adjustments that need to be made before this model can be used to achieve expected positive profit margins. Our results show that the only way to win GPP tournaments is to have a lineup that dramatically exceeds its projected fantasy score. As a result, individual lineups in these tournaments must have high variance. This means that for the probability distribution of fantasy points scored by any lineup, we want to increase the space under the tails. By increasing variance, we increase the probability that the lineup's actual score will be on either of the two tails of the lineup's fantasy point probability distribution. Practically, this means that the probability our lineup will do terribly *or* dramatically exceed expectations will increase. Consequently, increasing each lineup's variance in a contestant's portfolio increases the probability that one of the lineups will dramatically exceed expectations and win the tournament.

Methods to increase variance include adding constraints to the linear program model that incentivize creating lineups with players from the same team. This will increase variance because fantasy points overlap between hitters. For example, if Adam Eaton and Trea Turner are on base and Anthony Rendon hits a homerun and you have all three players on your team, your team will receive points for Adam Eaton scoring, for Trea Turner scoring, and for Anthony Rendon hitting a homerun. An individual lineup is much more likely to win a GPP tournament if a majority of the players in the lineup exceed expectations on that day. Players on the same team are more likely to exceed expectations together because of the fantasy point overlap.

Another possible method to increase individual lineup variance is to use data that projects what a player would score if they were to dramatically exceed expectations. For example, we can take data that projects the fantasy point score that a player will only exceed 10% of the time and optimize lineups based on those scores (e.g. if Mike Trout is expected to score 35 points or more 10% of the time, the value of 35 would be used in the model instead of the point value that we expect Mike Trout to score). The value of using this data instead of traditional projections is that it places an emphasis on players that are the most likely to significantly outperform their expectations. Traditionally, fantasy industry analysts refer to these types of players as "boom or bust" because they will either score exceedingly well and "boom" or score exceedingly poorly and "bust." Players that consistently score around their average are less likely to "boom" and thus should be less prioritized for GPP tournaments.

Although this model currently has its limitations, there are clear next steps that can be taken to substantially boost the model's ability to provide expected positive profit margins when competing in MLB GPP tournaments. Daily fantasy baseball has become a billion dollar industry with many websites offering paid services to help contestants optimize their lineups. Any model that provides consistent success in winning daily fantasy tournaments would be highly valuable in this industry. There are also clear parallels between daily fantasy lineup portfolios and stock portfolios or portfolios for mergers and acquisitions. If a linear programming approach that increases individual asset variance while also prioritizing overall asset diversity is successful in daily fantasy sports, there is evidence to suggest that pursuing such an approach with stocks or mergers and acquisitions could also be successful, given the similarities between daily fantasy and these market systems.